\documentclass{article}
\usepackage[utf8]{inputenc}
\usepackage{amsfonts,algorithmic,mathtools}
\usepackage{graphicx,overpic,epstopdf,xcolor}
\usepackage{cite}
\usepackage[section]{placeins}
\usepackage{amsthm,amsmath,amssymb}
\usepackage{hyperref,enumitem}
\usepackage{indentfirst,cite}
\ifpdf
  \DeclareGraphicsExtensions{.eps,.pdf,.png,.jpg}
\else
  \DeclareGraphicsExtensions{.eps}
\fi

\def\bb{\mathbb}
\def\cal{\mathcal}

\newtheorem{theorem}{Theorem}[section]

\theoremstyle{definition}
\newtheorem{definition}[theorem]{Definition}

\title{Estimating the Long-term Behavior of Biologically Inspired Agent-based Models}
\author{Daniel A. Cruz\thanks{Department of Medicine, University of Florida, Gainesville, FL. (\href{mailto:daniel.cruz@medicine.ufl.edu}{daniel.cruz@medicine.ufl.edu})}
\and Jack Toppen \thanks{School of Biological Sciences, Georgia Institute of Technology, Atlanta, GA. (\href{mailto:jtoppen3@gatech.edu}{jtoppen3@gatech.edu}, \href{mailto:epark90@gatech.edu}{epark90@gatech.edu})}
\and Eunbi Park \footnotemark[3]
\and Melissa L. Kemp \thanks{Wallace H. Coulter Department of Biomedical Engineering, Georgia Institute of Technology and Emory University, Atlanta, GA. (\href{mailto:melissa.kemp@bme.gatech.edu}{melissa.kemp@bme.gatech.edu})}
\and Elena Dimitrova \thanks{Mathematics Department, California Polytechnic State University, San Luis Obispo, CA. (\href{mailto:edimitro@calpoly.edu}{edimitro@calpoly.edu})}}
\date{November 30, 2022}

\begin{document}

\maketitle

\begin{abstract}
An agent-based model (ABM) is a computational model in which the local interactions of autonomous agents with each other and with their environment give rise to global properties within a given domain.
As the detail and complexity of these models has grown, so too has the computational expense of running several simulations to perform sensitivity analysis and evaluate long-term model behavior.
Here, we generalize a framework for mathematically formalizing ABMs to explicitly incorporate features commonly found in biological systems: appearance of agents (birth), removal of agents (death), and locally dependent state changes.
We then use our broader framework to extend an approach for estimating long-term behavior without simulations, specifically changes in population densities over time.
The approach is probabilistic and relies on treating the discrete, incremental update of an ABM via ``time steps'' as a Markov process to generate expected values for agents at each time step.
As case studies, we apply our extensions to both a simple ABM based on the Game of Life and a published ABM of rib development in vertebrates.
\end{abstract}

\textbf{Keywords:} agent-based model, cellular automata, population dynamics, global recurrence rule, interaction neighborhood

\textbf{MSC Codes:} 03D20, 60J05, 68Q80, 68U20, 92B05, 92C15

\section{Introduction}\label{s:intro}
In the broadest sense, an {\em agent-based model} (ABM) is a computational model of a target system in which autonomous agents governed by local rules interact with one another and with their environment \cite{Laubenbacher2012,Weimer2016}.
These local rules determine how an agent's current state updates at discrete, incremental ``time steps'' over the course of an ABM simulation.
The objective of ABM development is typically to connect local agent interactions with system-level properties across both time and space such as population dynamics and self-organization.
ABMs are found in virtually all areas of research \cite{Hinkelmann2011,Laubenbacher2012,Weimer2016}, supported by a growing list of software \cite{Abar2017} that range from general-purpose platforms such as NetLogo \cite{Wilensky1999} and the PythonABM library \cite{Toppen2022} to specialized platforms like CompuCell3D \cite{Swat2012}, Morpheus \cite{Starruss2014}, and PhysiCell \cite{Ghaffarizadeh2018} which aim to capture as many system details as possible in their implementation.
In particular, the use of ABMs in the context of cellular biology has grown rapidly over the last decade \cite{Metzcar2019,Glen2019,Buttenshon2020,Wauford2022} because of how cells can naturally be considered autonomous agents.
The agent-based modeling of cellular interactions, specifically of morphogenesis \cite{Glen2019}, is of interest because of how ABMs can incorporate processes such as cell division (e.g.  mitosis) and death (e.g. apoptosis or necrosis) while incorporating the spatial information of each cell that most other systems biology models cannot \cite{Machado2011,Bartocci2016,Cruz2021}.

As ABMs become more ubiquitous, there is a growing interest in the development of methods to analyze these models.
In most cases, the evaluation of ABMs against their associated systems is either too qualitative or too computationally expensive to be suitable for research fields like systems biology where the throughput of data acquisition is increasing while costs are decreasing.
Qualitative evaluation typically relies on a visual comparison between simulated and collected data while quantitative evaluation involves running several simulations and performing statistical analysis to generate a measure of model ``consistency'' or ``reproducibility''; the latter often requires significant computational resources and long simulation run-times \cite{Laubenbacher2012,Nardini2021a}.
As an alternative, several mathematical frameworks have been proposed over the years for formalizing ABMs, including finite dynamical systems, cellular automata, and Markov chains \cite{Hinkelmann2011,Laubenbacher2012,Yereniuk2019,KhudaBukhsh2019,Nardini2021a}.
The common goal of these frameworks is to estimate the long-term behavior of an ABM (e.g. the change in agent populations over time) without needing data from simulations.
However, these frameworks (i) fix the number of agents across a simulation and/or (ii) ignore or discretize the positional information of agents.
Because of the aforementioned restrictions, such frameworks cannot easily be applied to ABMs which (i) simulate biological phenomena involving agent birth and death and (ii) exhibit emergent spatial properties such as pattern formation.

In this work, we extend a framework that (i) explicitly allows an agent's location to be continuous or discrete and (ii) provides a clear connection between this positional information and the expected population densities of each (agent) state \cite{Yereniuk2019}.
This framework requires the update process by which local rules modify agents at each time step to be Markovian in order to simplify later computations.
We note that this condition is relatively mild and agrees with the way that most ABMs are coded in practice.
Our extension of this framework allows for the number of agents to change over time according to an ABM's local rules; as such, our definition of an ABM accommodates a broad set of models actively used in research.
We also generalize the computation of a state's ``global recurrence rule'' (GRR) \cite{Yereniuk2019} to calculate changes in population density for that state based on the ABM's local rules and an initial set of model parameters.
We present an approach for estimating a state's GRR rather than explicitly calculating it to focus on practical application and provide two case studies of our approach: an ABM based on the Game of Life cellular automaton \cite{Gardner1970} (which we call ``GoL-like'') and an ABM developed with experimental data to study early rib development in vertebrates \cite{Fogel2017}.
Our approach relies on both describing an ABM's local rules in terms of an agent's neighborhood and using agents' neighborhoods to approximate expected behavior across the ABM's environment.

In Section~\ref{s:def}, we present the definitions and notation for our generalized framework.
In Section~\ref{s:GoL}, we focus on a GoL-like ABM which serves as an illustrative example of both the framework and our approach to approximating a state's GRR.
In Section~\ref{s:ribABM}, we then apply our approach to a published ABM for rib development and provide estimates of population changes for three agent states (i.e. cell types) across four experimental trials.
Both ABMs were implemented in NetLogo \cite{Wilensky1999} and their associated GRR computations were written and run in MATLAB \cite{matlab2022}.
See \url{https://github.com/kemplab/ABM-Math-Framework} for all of the code associated to this paper.

\section{Definitions}\label{s:def}

For any set $S$, we use $\bb{M}(S)$ to denote the set of all finite multisets (or collections) whose elements are in $S$.
For an element $a\in S$ and a finite multiset $M\in\bb{M}(S)$ containing $a$, we take $M\setminus a$ to mean the multiset resulting from the removal of one copy of $a$ from $M$, keeping in mind that there may be several copies in $M$.
Though a finite multiset can always be written as a set by using additional notation to distinguish duplicate elements if necessary, we prefer to refer to such objects as multisets for simplicity and use square brackets for multisets accordingly.

For an interval $I\subset\bb{R}$, we use $Uniform(I)$ to denote the uniform probability distribution over $I$.
While we also use square brackets for some intervals in $\bb{R}$, the context and notation below make it easy to distinguish between multisets and intervals.
For $n\in\bb{Z}^+$, $\bb{R}^n$ is the classic $n$-dimensional Euclidean space with a fixed origin $O_n$, typically denoted as $O$.

\subsection{Formalizing Agent-based Models}\label{ss:ABM}
We primarily focus on extending the definitions of \cite{Yereniuk2019} in the context of when the ABM environment (i.e. bounded region of interest) $\Omega$ is a connected, bounded subset of $\bb{R}^n$.
If $\Omega$ is instead a graph $G=(V,E)$, then we can employ standard graph theory definitions and adjust our discussion accordingly.
We begin by formally defining an agent below, choosing to represent it as a point in the environment $\Omega$.
Note that we can easily extend our definitions to represent an agent as having a ``shape'' if we wish; see Appendix~\ref{a:def} for details.

\begin{definition}\label{def:agent}
Let $\Sigma$ be a finite set of {\em states}, and let $\Omega$ be (1) a connected, bounded subset of $\bb{R}^n$ or (2) a graph $(V,E)$.
An {\em agent} $\alpha=(s,p,\cal{N})$ is an ordered triple where $s\in\Sigma$, and $p$ and $\cal{N}$ are defined as follows with respect to $\Omega$:
\begin{itemize}
    \item $p$ is (1) a point in the set $\Omega$ or (2) a vertex in the graph $\Omega=(V,E)$, as appropriate.
    \item $\cal{N}$ is (1) a connected subset such that $p\in\cal{N}\subset\Omega$ or (2) $\cal{N}$ is a finite subset such that $p\in\cal{N}\subseteq V$. 
\end{itemize}
In either of the cases (1) or (2), we call $p$ the {\em position of $\alpha$}, $s$ the {\em state of $\alpha$}, and $\cal{N}$ the {\em neighborhood of $\alpha$}.
We use the notation $p(\alpha)$ to refer to $p$ and use the similar notation for $s$ and $\cal{N}$.
Finally, we say that $\Omega$ is an {\em environment} in this context and use $\Lambda(\Sigma,\Omega)$ to denote the set of all agents over $\Sigma$ and $\Omega$. 
\end{definition}

Given a multiset of agents $\cal{X}\in\bb{M}(\Lambda(\Sigma,\Omega))$ within the environment $\Omega$, we can now formalize the process by which an agent in $\cal{X}$ will update its attributes (i.e. position, state, and neighborhood) based on those of its ``neighbors'' in $\cal{X}$ and potentially produce new agents. 
Let $\Lambda(\Sigma,\Omega)$ be given for an environment $\Omega$ and a set of states of $\Sigma$.
\begin{itemize}
\item A {\em local transition rule (over $\Lambda(\Sigma,\Omega)$)} is a mapping $$f:\Lambda(\Sigma,\Omega)\times \bb{M}(\Lambda(\Sigma,\Omega))\to \Lambda(\Sigma,\Omega).$$
\item A {\em local production rule (over $\Lambda(\Sigma,\Omega)$)} is a mapping $$g:\Lambda(\Sigma,\Omega)\times \bb{M}(\Lambda(\Sigma,\Omega))\to \bb{M}(\Lambda(\Sigma,\Omega)).$$
\end{itemize}

\begin{figure}[ht]
\centering
\begin{overpic}[width=\textwidth]{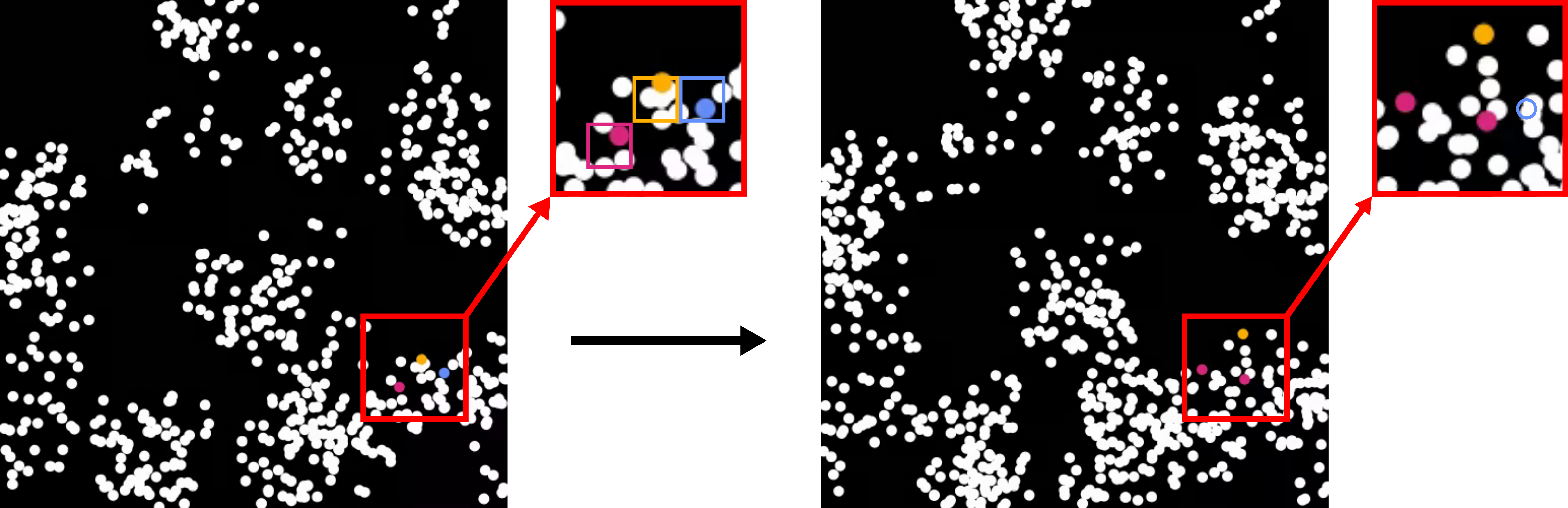}
\put(14,-3){time $t$}
\put(63, -3){time $t+1$}
\put(38,7){$f_\cal{G}$, $g_\cal{G}$}
\end{overpic}
\caption{A visualization of the local transition rule $f_\cal{G}$ and the local production rule $g_\cal{G}$ introduced in Section~\ref{ss:ABM} being applied to a multiset $\cal{X}\in\bb{M}(\Lambda(\Sigma_\cal{G},\Omega_\cal{G}))$ (left) whose agents all have neighborhoods of the form $[i,i+1)\times[j,j+1)$ for some integers $0\leq i,j\leq 19$.
We use unit circles to visualize agents in this context and let the larger black squares in which these agents are located be the environment $\Omega_\cal{G}$.
We take the position of each agent as the center of its corresponding unit circle.
To illustrate the dynamics of $f_\cal{G}$ and $g_\cal{G}$, we color all but three living agents white within the environment on the left and assume that there are no other agents.
The remaining three agents are colored magenta, orange, and blue; we also highlight their neighborhoods similarly within the inset on the left.
The blue agent only has one neighbor (left), so it will ``die'' after $f_\cal{G}$ is applied, visualized by making the agent transparent with only a blue outline (right).
On the other hand, the orange agent will survive and move (right) because it has five neighbors (left).
Note neither of these agents have the correct number of neighbors to produce offspring according to the definition of $g_\cal{G}$.
The magenta agent has 3 neighbors (left), so it will survive, move, and produce one identical offspring which itself will move; we color this offspring magenta as well (right).
We can consider this visualization as applying $f_\cal{G}$ and $g_\cal{G}$ at some time $t$ (left) during a simulation of the associated GoL-like ABM to produce an updated multiset at time $t+1$ (right); see Definition~\ref{def:ABM}.
}
\label{fig:GoL_rules}
\end{figure}

As a simple example, we describe an ABM based on the Game of Life cellular automaton \cite{Gardner1970}; see Figure~\ref{fig:GoL_rules} for a visual reference.
For this ABM, we have an environment $\Omega_\cal{G} = [0,20)\times[0,20) \subset \bb{R}^2$ and a set of states $\Sigma_\cal{G} = \{0,1\}$, where $0$ indicates ``dead'' and $1$ indicates ``alive''.
We define the local transition rule $f_\cal{G}$ and the local production rule $g_\cal{G}$ for this ABM as follows for an agent $\alpha=(s,p,\cal{N}) \in \Lambda(\Sigma_\cal{G},\Omega_\cal{G})$ and a multiset $\cal{X}\in\bb{M}(\Lambda(\Sigma_\cal{G},\Omega_\cal{G}))$:
\begin{itemize}
    \item $f_\cal{G}(\alpha,\cal{X}) = \alpha$ if $\alpha\not\in\cal{X}$ or if there exists an agent $\alpha'$ in $\cal{X}$ such that $\cal{N}(\alpha')\neq [i,i+1)\times[j,j+1)$ for some integers $0\leq i,j\leq 19$.
    Otherwise,
    $$f_\cal{G}(\alpha,\cal{X})=\begin{cases}
    (1,p',\cal{N}') & \text{if } s(\alpha)=1 \text{ and }\\
    & 2\leq\big\vert[\beta\in\cal{X}\setminus\alpha \mid s(\beta)=1 \text{ \& } p(\beta)\in\cal{N}(\alpha) ]\big\vert\leq 8\\
    (0,p,\cal{N}) & \text{otherwise,}
    \end{cases}$$
    where 
    $$
    p'=\begin{cases}
    p+(\cos(\theta),\sin(\theta)) & \text{if } p+(\cos(\theta),\sin(\theta)) \in \Omega_\cal{G}\\
    p & \text{otherwise}
    \end{cases}
    $$
    for $\theta\sim Uniform[0,2\pi)$ and $\cal{N}' = [i,i+1)\times[j,j+1)$ for integers $0\leq i,j\leq 19$ such that $p'\in\cal{N}'$ (as required).
    \item $g_\cal{G}(\alpha,\cal{X}) = \emptyset$ if $\alpha\not\in\cal{X}$ or if there exists an agent $\alpha'$ in $\cal{X}$ such that $\cal{N}(\alpha')\neq [i,i+1)\times[j,j+1)$ for some integers $0\leq i,j\leq 19$.
    Otherwise,
    $$g_\cal{G}(\alpha,\cal{X})=\begin{cases}
    \{(1,p',\cal{N}')\} & \text{if } s(\alpha)=1 \text{ and }\\
    & 2\leq\big\vert[\beta\in\cal{X}\setminus\alpha \mid s(\beta)=1 \text{ \& } p(\beta)\in\cal{N}(\alpha) ]\big\vert\leq 4\\
    \emptyset & \text{otherwise,}
    \end{cases}$$
    where $p'$ and $\cal{N}'$ are defined as in $f_\cal{G}$.
\end{itemize}

Note that the common restriction in the definitions above is that $f_\cal{G}$ and $g_\cal{G}$ do ``nothing'' unless $\alpha\in\cal{X}$ and the neighborhood of each agent in $\cal{X}$ is (essentially) a unit square\footnote{We refer to sets of the form $[i,i+1)\times[j,j+1)$ as {\em unit squares} in this work for simplicity, noting that this deviation from standard terminology does not affect our results in Section~\ref{s:GoL}.}.
Informally, we can describe the local transition and production rules of this ABM as follows for an agent $\alpha$ in $\Omega_\cal{G}$ and a (valid) multiset $\cal{X}$ containing $\alpha$:
\begin{itemize}
    \item $\alpha$ only survives if it has 2 to 8 (living) neighbors, in which case $\alpha$ moves a distance of at most 1 within $\Omega_\cal{G}$ in a random direction.
    \item $\alpha$ only produces a (living) offspring if it has 2 to 4 (living) neighbors, in which case its offspring moves at random as above.
\end{itemize}
In Figure~\ref{fig:GoL_rules}, we visualize the application of $f_\cal{G}$ and $g_\cal{G}$ on each agent in a multiset $\cal{X}\in\bb{M}(\Lambda(\Sigma_\cal{G},\Omega_\cal{G}))$ whose neighborhoods are all unit squares.

In general, given a local transition rule $f$ and any multiset $\cal{X}\in\bb{M}(\Lambda(\Sigma,\Omega))$, we use $f(\cal{X},\cal{X})$ to denote the multiset $\cal{X}'=[f(\alpha,\cal{X})\mid \alpha\in\cal{X}]$; informally, $\cal{X}'$ is the multiset of agents which we obtain from applying the local transition rule to every agent $\alpha$ in $\cal{X}$.
On the other hand, a local production rule $g$ yields a multiset of agents $g(\cal{X},\cal{X})=\bigcup_{\alpha\in\cal{X}} g(\alpha,\cal{X})$ which will be ``added'' to the multiset $f(\cal{X},\cal{X})$, as described in the next definition.

\begin{definition}\label{def:ABM}
Let $\Lambda(\Sigma,\Omega)$ be given for an environment $\Omega$ and a set of states $\Sigma$.
An {\em agent-based model} (ABM) $\cal{A}=(\Sigma,\Omega,f,g)$ is a 4-tuple where $f$ and $g$ are local transition and production rules, respectively, over $\Lambda(\Sigma,\Omega)$.
We collectively refer to $f$ and $g$ as the {\em update rules} of $\cal{A}$.

Given $\cal{A}$ and a multiset $\cal{X}_0\in\bb{M}(\Lambda(\Sigma,\Omega))$, we define the sequence of (finite) multisets $\{\cal{X}_t\}_{t\geq 0}$ as follows for $t\geq 1$:
\begin{equation*}
\cal{X}_t = f(\cal{X}_{t-1},\cal{X}_{t-1})\cup g(\cal{X}_{t-1},\cal{X}_{t-1}) \in \bb{M}(\Lambda(\Sigma,\Omega)).
\end{equation*}
Informally, each $\cal{X}_t$ is obtained from $\cal{X}_{t-1}$ by joining the multiset $f(\cal{X}_{t-1},\cal{X}_{t-1})$ yielded by the transition rule $f$ with the multiset $g(\cal{X}_{t-1},\cal{X}_{t-1})$ yielded by the production rule $g$.
We say that the sequence $\{\cal{X}_t\}_{t\geq 0}$ is a {\em simulation (of $\cal{A}$)} and refer to $\cal{X}_0$ as an {\em initialization (of $\cal{A}$)} in this context.
\end{definition}

In the definition above, a simulation $\{\cal{X}_t\}_{t\geq 0}$ directly corresponds to a computational simulation in the standard context of discussing and developing ABMs.
It is common to refer to an element $\cal{X}_i$ of $\{\cal{X}_t\}_{t\geq 0}$ as ``the simulation at time (step) $t=i$''; thus, we adopt the convention of referring to the subscript $i$ as the {\em time (step)} in the context of simulations.
Note that we may now formally define the GoL-like ABM $\cal{A}_\cal{G}$ described in Figure~\ref{fig:GoL_rules}: $\cal{A}_\cal{G}=(\Sigma_\cal{G},\Omega_\cal{G},f_\cal{G},g_\cal{G})$ where $\Omega_\cal{G} = [0,20)\times[0,20) \subset \bb{R}^2$, $\Sigma_\cal{G} = \{0,1\}$, and the local update rules $f_\cal{G}$ and $g_\cal{G}$ are as defined earlier in this section.
With this definition, we can formally consider the two multisets in Figure~\ref{fig:GoL_rules} as being consecutive elements of a simulation of $\cal{A}_\cal{G}$.

Before getting into definitions associated with analyzing long-term behavior of ABMs, we note that our definition of a local transition rule $f$ differs from the one in \cite{Yereniuk2019} because it is deterministic instead of stochastic.
For brevity, we modify our definitions here to simplify the construction and analysis of our first example in Section~\ref{s:GoL}.
Nonetheless, it is straightforward to define both local transition and production rules as stochastic mappings and to subsequently change all associated definitions thereafter; see Appendix~\ref{a:def} for details.

\subsection{Global Recurrence Rules for ABMs}\label{ss:GRR}
We now turn our attention to calculating state density changes during the simulation of a given ABM using ``transition'' and ``production'' regions.

\begin{definition}\label{def:transition_region}
Let $\cal{A}=(\Sigma,\Omega,f,g)$ be an ABM and $\{\cal{X}_t\}_{t\geq 0}$ be a simulation of $\cal{A}$.
Let $t\geq 0$, $x\in\Omega$, and $\cal{V},\cal{U}$ be states in $\Sigma$.
We say that $x$ is a {\em $(\cal{V},\cal{U})$-transition point at time $t$} if for any agent $\alpha\in \Lambda(\Sigma,\Omega)$ with $s(\alpha)=\cal{V}$ and $p(\alpha)=x$, we have that $s\left(f(\alpha,\cal{X}_t)\right)=\cal{U}$.
The {\em $(\cal{V},\cal{U})$-transition region at time $t$} is the set 
\begin{equation*}
B_t^{\cal{V},\cal{U}}=\{x\in\Omega\mid x \text{ is a } (\cal{V},\cal{U}) \text{-transition point at time } t\}.
\end{equation*}
\end{definition}

Note that the set of states $\Sigma$ for an ABM $\cal{A}=(\Sigma,\Omega,f,g)$ can have an element $\epsilon$ designated as the {\em death state} if agents are allowed to ``die'' or otherwise disappear, as with the GoL-like ABM $\cal{A}_\cal{G}$ introduced in Section~\ref{ss:ABM}.
In this case, $f$ and $g$ are defined such that $f(\alpha,\cal{X})=\alpha$ and $g(\alpha,\cal{X})=\emptyset$ whenever $s(\alpha)=\epsilon$ (i.e. whenever the agent is ``dead'').
The $(\epsilon,\cal{U})$-transition region at time $t$ is as follows: $B_t^{\epsilon,\epsilon}=\Omega$ (i.e. when $\cal{U}=\epsilon$) and $B_t^{\epsilon,\cal{U}}=\emptyset$ otherwise.
Informally, this means that a ``dead'' agent must remain in this state for the rest of the simulation $\{\cal{X}_t\}_{t\geq 0}$.

\begin{definition}\label{def:production_region}
Let $\cal{A}=(\Sigma,\Omega,f,g)$ be an ABM and $\{\cal{X}_t\}_{t\geq 0}$ be a simulation of $\cal{A}$.
Let $t\geq 0$, $x\in\Omega$, and $\cal{V},\cal{U}$ be states in $\Sigma$.
We say that $x$ is a {\em $(\cal{V},\cal{U})$-production point at time $t$} if for any agent $\alpha\in \Lambda(\Sigma,\Omega)$ with $s(\alpha)=\cal{V}$ and $p(\alpha)=x$, we have that $\exists \beta\in g(\alpha,\cal{X}_t)$ with $s(\beta)=\cal{U}$.
The {\em $(\cal{V},\cal{U})$-production region at time $t$} is the set 
\begin{equation*}
C_t^{\cal{V},\cal{U}}=\{x\in\Omega\mid x \text{ is a } (\cal{V},\cal{U}) \text{-production point at time } t\}.
\end{equation*}
\end{definition}

Given an ABM $\cal{A}=(\Sigma,\Omega,f,g)$ and a simulation $\{\cal{X}_t\}_{t\geq 0}$ of $\cal{A}$, we can repeat the following observation from \cite{Yereniuk2019} for any $t\in\bb{Z}^+$, $\alpha\in\cal{X}_t$, and $\cal{V},\cal{U}\in\Sigma$:
\begin{equation}
\bb{P}\left(s\left(f(\alpha,\cal{X}_t)\right)=\cal{U} \;\middle\vert\; s(\alpha)=\cal{V}\right) = \bb{P}\left(p(\alpha)\in B_t^{\cal{V},\cal{U}}\right).
\end{equation}\label{eqn:YO_obs1}
We can make a similar observation regarding the local production rule for any $\beta\in g(\alpha,\cal{X}_t)$:
\begin{equation}
\bb{P}\Big(s(\beta)=\cal{U} \;\Big\vert\; s(\alpha)=\cal{V} \text{ and } \beta\in g(\alpha,\cal{X}_t)\Big) = \bb{P}\left(p(\alpha)\in C_t^{\cal{V},\cal{U}}\right).
\end{equation}\label{eqn:YO_obs2}

We define the mapping $\cal{D}_\cal{U}:\bb{Z}^{\geq 0}\to\bb{Z}^{\geq 0}$ for each $\cal{U}\in\Sigma$ such that $\cal{D}_\cal{U}(t)=\left\lvert\left[\alpha\in\cal{X}_t\;\middle\vert\;s(\alpha)=\cal{U}\right]\right\rvert$.
Then by adapting the work from \cite{Yereniuk2019} to include agent production, we can determine that the expected number of agents with state $\cal{U}$ at time $(t+1)\in\bb{Z}^+$, denoted $\bb{E}\left(\cal{D}_\cal{U}(t+1)\right)$:

\begin{align*}
&\bb{E}\left(\cal{D}_\cal{U}(t+1)\right) = \bb{E}\left(\left\lvert\left[\alpha\in\cal{X}_{t+1}\;\middle\vert\;s(\alpha)=\cal{U}\right]\right\rvert\right)\\
&= \sum_{\alpha\in\cal{X}_t}\left[\bb{P}\Big(s\left(f(\alpha,\cal{X}_t)\right)=\cal{U}\Big) + \quad\smashoperator{\sum_{\beta\in g(\alpha,\cal{X}_t) \text{ s.t. } s(\beta)=\cal{U}}}\quad \bb{P}\Big(s(\beta)=\cal{U} \;\Big\vert\; \beta\in g(\alpha,\cal{X}_t) \Big)\right]\\
&= \sum_{\cal{V}\in\Sigma}\sum_{[\alpha\in\cal{X}_t\vert(\alpha)=\cal{V}]}\bigg(\bb{P}\Big(s\left(f(\alpha,\cal{X}_t)\right)=\cal{U}\;\Big\vert\;s(\alpha)=\cal{V}\Big) +\\
&\qquad\qquad\qquad\Big\vert[\beta\in g(\alpha,\cal{X}_t)\;\vert\; s(\beta)=\cal{U}]\Big\vert\cdot\bb{P}\Big(s(\beta)=\cal{U} \;\Big\vert\; s(\alpha)=\cal{V} \text{ and } \beta\in g(\alpha,\cal{X}_t) \Big)\bigg)\\
&= \sum_{\cal{V}\in\Sigma}\sum_{[\alpha\in\cal{X}_t\vert(\alpha)=\cal{V}]} \bb{P}\left(p(\alpha)\in B_t^{\cal{V},\cal{U}}\right) + \Big\vert[\beta\in g(\alpha,\cal{X}_t)\;\vert\; s(\beta)=\cal{U}]\Big\vert\cdot\bb{P}\left(p(\alpha)\in C_t^{\cal{V},\cal{U}}\right)
\end{align*}

The equations above provide us with the main method of calculating long-term ABM behavior using this framework, which we summarize in the following definition.

\begin{definition}\label{def:GRR}
Let $\cal{A}=(\Sigma,\Omega,f,g)$ be an ABM and $\{\cal{X}_t\}_{t\geq 0}$ be a simulation of $\cal{A}$.
For $\cal{U}\in\Sigma$ and $t\geq 0$, the {\em global recurrence rule (GRR) of $\cal{U}$ (with respect to $t+1$)}, $\bb{E}\left(\cal{D}_\cal{U}(t+1)\right)$, is given by the following expression:
\begin{equation}
\sum_{\cal{V}\in\Sigma}\smashoperator[r]{\sum_{[\alpha\in\cal{X}_t\vert(\alpha)=\cal{V}]}}\quad \bb{P}\left(p(\alpha)\in B_t^{\cal{V},\cal{U}}\right) + \Big\vert[\beta\in g(\alpha,\cal{X}_t)\;\vert\; s(\beta)=\cal{U}]\Big\vert\cdot\bb{P}\left(p(\alpha)\in C_t^{\cal{V},\cal{U}}\right).
\end{equation}
\end{definition}

By Definition~\ref{def:GRR}, note that finding the GRR of a state $\cal{U}$ comes down to determining $\bb{P}\left(p(\alpha)\in B_t^{\cal{V},\cal{U}}\right)$ and $\bb{P}\left(p(\alpha)\in C_t^{\cal{V},\cal{U}}\right)$ for some $\alpha\in\cal{X}_t$.
It follows that using the GRR to calculate or estimate long-term ABM behavior works best in instances when these probabilities can be determined or approximated, respectively.
We provide examples of such ABMs in the Sections \ref{s:GoL} and \ref{s:ribABM}.

\section{Game of Life-like ABM}\label{s:GoL}

\begin{figure}[ht]

\begin{minipage}{0.23\textwidth}
\centering
\begin{overpic}[width=\linewidth]{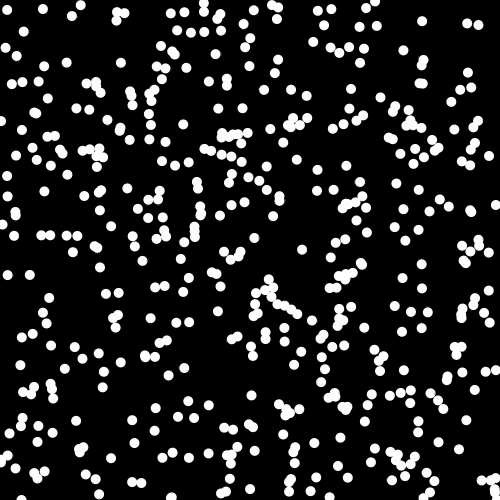}
\put(38,-13){$t=0$}
\end{overpic}

\vspace{1.5\baselineskip}

\begin{overpic}[width=\linewidth]{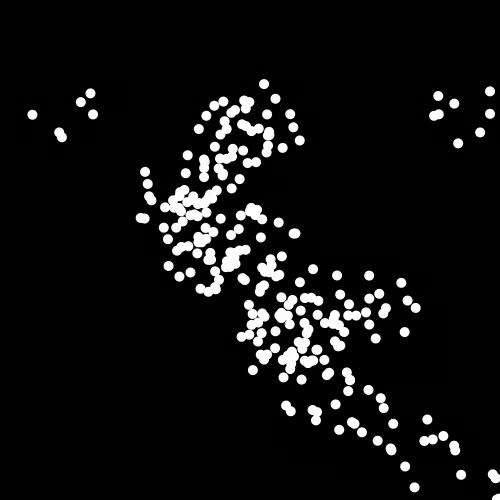}
\put(38,-13){$t=4$}
\end{overpic}

\end{minipage}
\hfill
\begin{minipage}{0.23\textwidth}
\centering
\begin{overpic}[width=\linewidth]{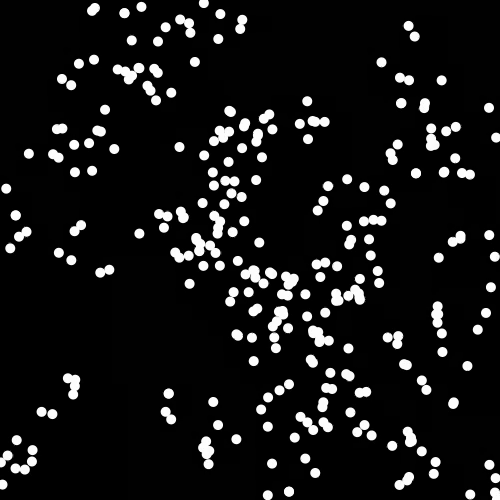}
\put(38,-13){$t=1$}
\end{overpic}

\vspace{1.5\baselineskip}

\begin{overpic}[width=\linewidth]{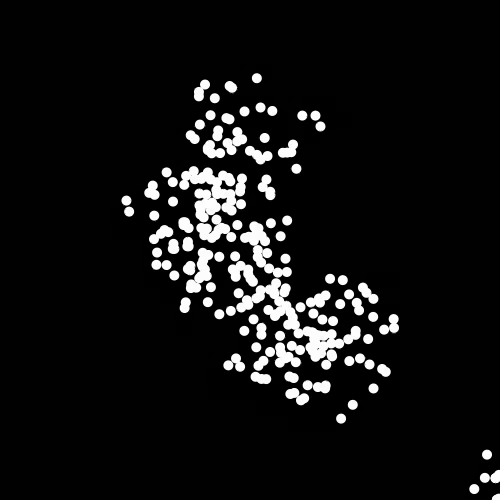}
\put(38,-13){$t=5$}
\end{overpic}

\end{minipage}
\hfill
\begin{minipage}{0.23\textwidth}
\centering
\begin{overpic}[width=\linewidth]{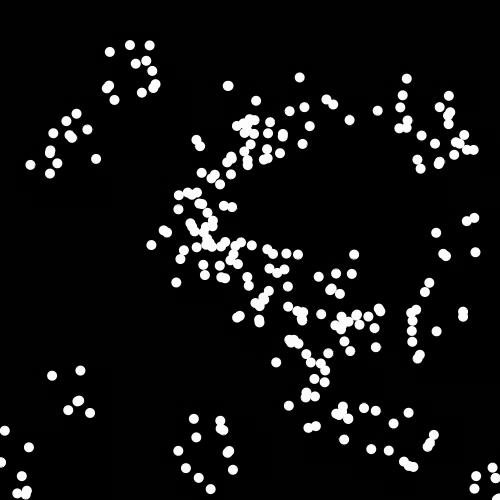}
\put(38,-13){$t=2$}
\end{overpic}

\vspace{1.5\baselineskip}

\begin{overpic}[width=\linewidth]{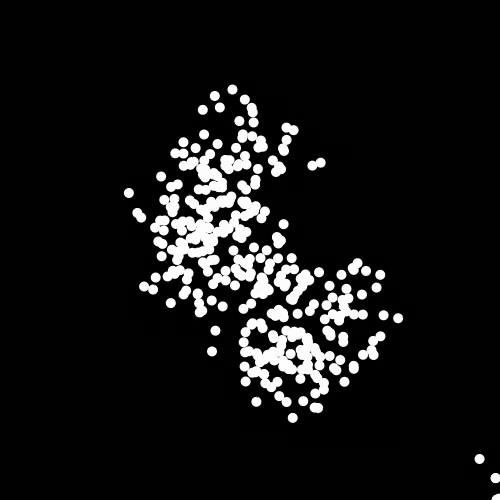}
\put(38,-13){$t=6$}
\end{overpic}

\end{minipage}
\hfill
\begin{minipage}{0.23\textwidth}
\centering
\begin{overpic}[width=\linewidth]{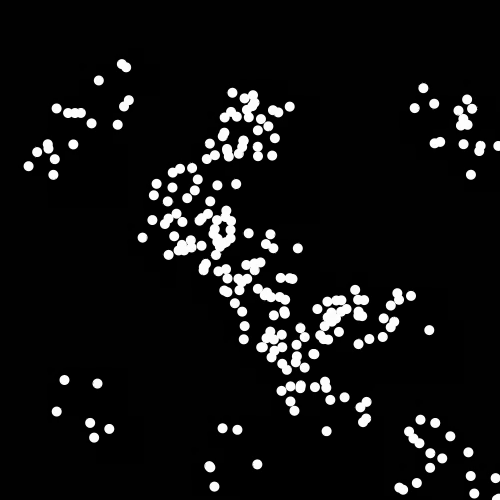}
\put(38,-13){$t=3$}
\end{overpic}

\vspace{1.5\baselineskip}

\begin{overpic}[width=\linewidth]{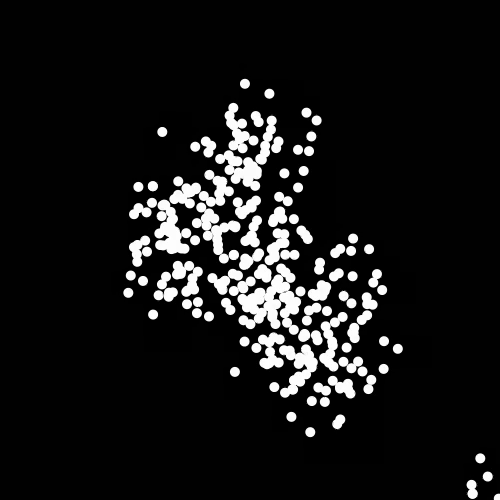}
\put(38,-13){$t=7$}
\end{overpic}

\end{minipage}

\vspace{1.5\baselineskip}

\caption{The first eight elements of a simulation $\{\cal{X}_t\}_{t\geq 0}$ of $\cal{A}_\cal{G}$ labeled according to their time step $t$.
For this simulation, there are $n_0=500$ agents in initialization $\cal{X}_0$ which are ``alive'' (i.e. have a state of 1) and are located uniformly at random in the environment $\Omega_\cal{G}$ of $\cal{A}_\cal{G}$.
No other agents exist in $\cal{X}_0$.
We use white unit circles to visualize living agents as in Figure~\ref{fig:GoL_rules} and do not visualize dead agents at all, as is standard for most ABMs.}
\label{fig:GoL_simulation}
\end{figure}

In this section, we consider the GoL-like ABM $\cal{A}_\cal{G}$ introduced in Section~\ref{ss:ABM} further as an example ABM with simple rules.
We will approximate the GRR of living agents within $\cal{A}_\cal{G}$ under the conditions that living agents in the initialization of a simulation $\{\cal{X}_t\}_{t\geq 0}$ are located in the environment uniformly at random.
We begin by establishing notation to describe a broader class of GoL-like ABMs based on $\cal{A}_\cal{G}$.

\begin{definition}\label{def:GoL}
Let $w, \ell_{surv}, u_{surv}, \ell_{rep}, u_{rep} \in\bb{Z}^{\geq 0}$ such that $\ell_{surv}\leq u_{surv}$ and $\ell_{rep}\leq u_{rep}$.
Suppose that $\cal{A}=(\Sigma,\Omega,f,g)$ is an ABM with states $\Sigma = \{0,1\}$, environment $\Omega=[0,w)\times[0,w)\subset\bb{R}^2$, and local update rules $f$ and $g$ defined as follows for an agent $\alpha=(s,p,\cal{N}) \in \Lambda(\Sigma,\Omega)$ and a multiset $\cal{X}\in\bb{M}(\Lambda(\Sigma,\Omega))$:
\begin{itemize}
    \item $f(\alpha,\cal{X}) = \alpha$ if $\alpha\not\in\cal{X}$ or if there exists an agent $\alpha'$ in $\cal{X}$ such that $\cal{N}(\alpha')\neq [i,i+1)\times[j,j+1)$ for some integers $0\leq i,j\leq w-1$.
    Otherwise,
    $$f(\alpha,\cal{X})=\begin{cases}
    (1,p',\cal{N}') & \text{if } s(\alpha)=1 \text{ and }
    \ell_{surv}\leq\big\vert[\beta\in\cal{X}\setminus\alpha \mid s(\beta)=1 \text{ \&}\\
    & \qquad\qquad\qquad\qquad\qquad\quad p(\beta)\in\cal{N}(\alpha) ]\big\vert\leq u_{surv}\\
    (0,p,\cal{N}) & \text{otherwise,}
    \end{cases}$$
    where 
    $$
    p'=\begin{cases}
    p+(\cos(\theta),\sin(\theta)) & \text{if } p+(\cos(\theta),\sin(\theta)) \in \Omega\\
    p & \text{otherwise}
    \end{cases}
    $$
    for $\theta\sim Uniform[0,2\pi)$ and $\cal{N}' = [i,i+1)\times[j,j+1)$ for integers $0\leq i,j\leq w-1$ such that $p'\in\cal{N}'$ (as required).
    \item $g(\alpha,\cal{X}) = \emptyset$ if $\alpha\not\in\cal{X}$ or if there exists an agent $\alpha'$ in $\cal{X}$ such that $\cal{N}(\alpha')\neq [i,i+1)\times[j,j+1)$ for some integers $0\leq i,j\leq w-1$.
    Otherwise,
    $$g(\alpha,\cal{X})=\begin{cases}
    \{(1,p',\cal{N}')\} & \text{if } s(\alpha)=1 \text{ and }
    \ell_{rep}\leq\big\vert[\beta\in\cal{X}\setminus\alpha \mid s(\beta)=1 \text{ \&}\\
    & \qquad\qquad\qquad\qquad\qquad\quad p(\beta)\in\cal{N}(\alpha) ]\big\vert\leq u_{rep}\\
    \emptyset & \text{otherwise,}
    \end{cases}$$
    where $p'$ and $\cal{N}'$ are defined as in $f$.
\end{itemize}
Then we say that $\cal{A}$ is a {\em Game of Life-like (GoL-like) ABM} and use $\cal{G}(w,$ $\ell_{surv},$ $u_{surv},$ $\ell_{rep},$ $u_{rep})$ to denote it.
\end{definition}

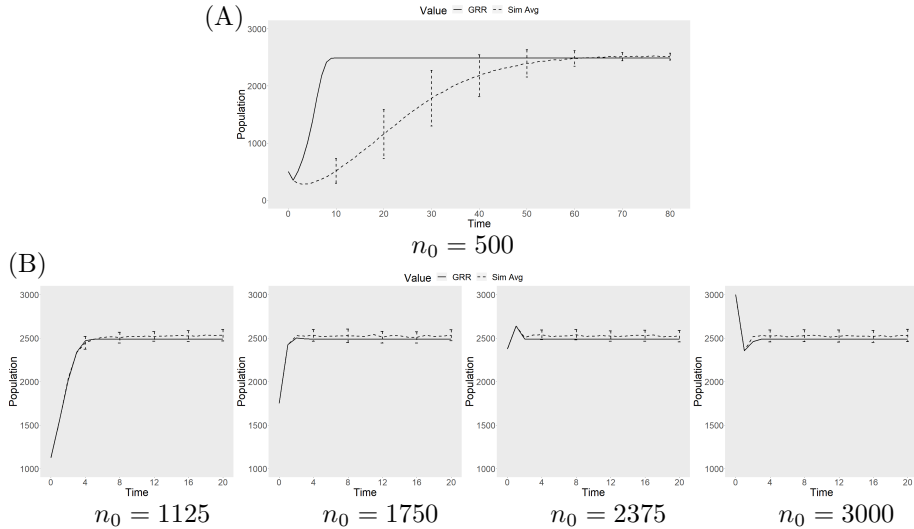
\begin{figure}[ht]
    \centering
    \begin{overpic}[width=0.5\textwidth]{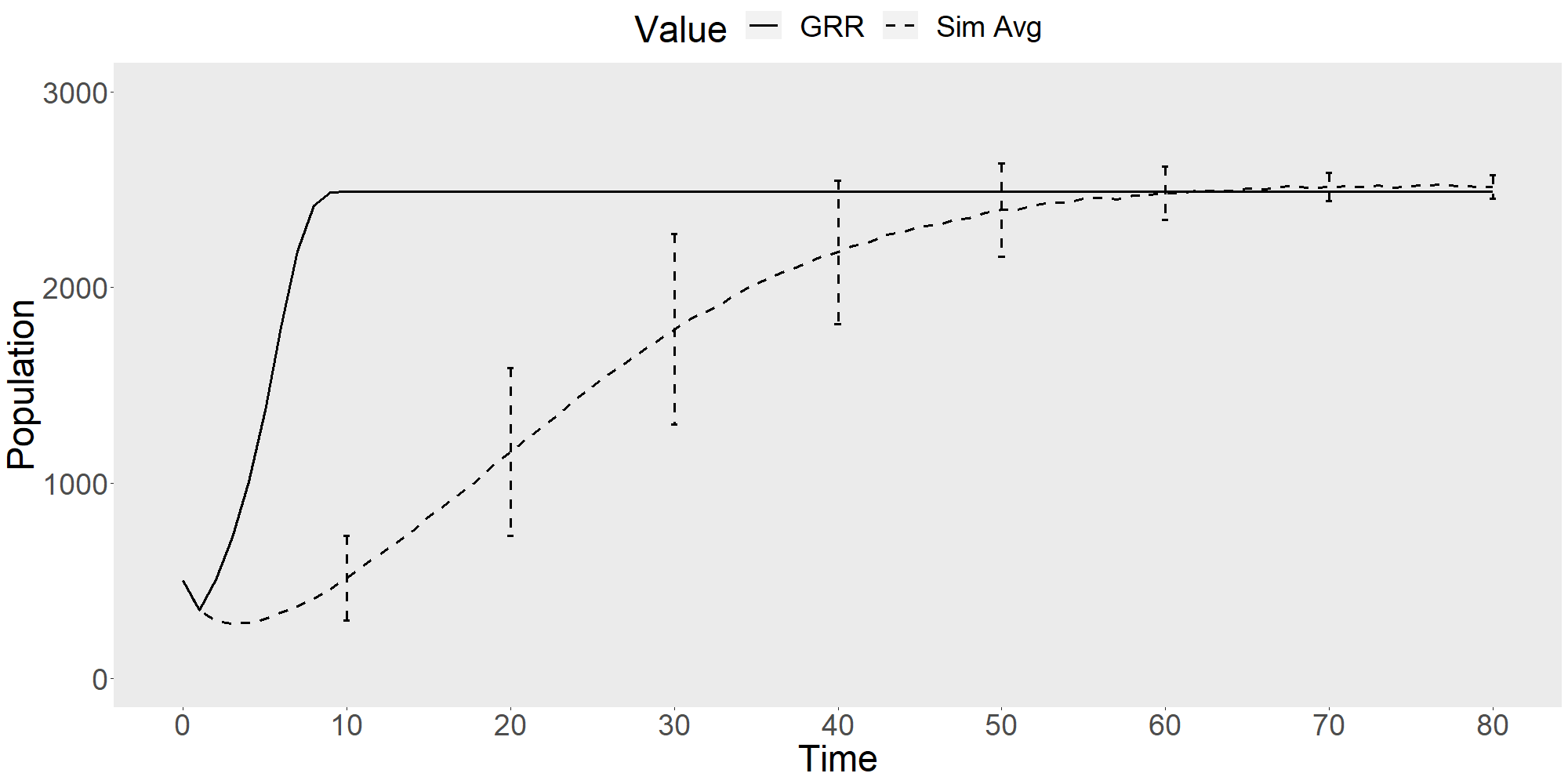}
    \put(-7,45){(A)}
    \put(38.5,-5){$n_0=500$}
    \end{overpic}
    
    \vspace{1.5\baselineskip}
    
    \begin{overpic}[width=\textwidth]{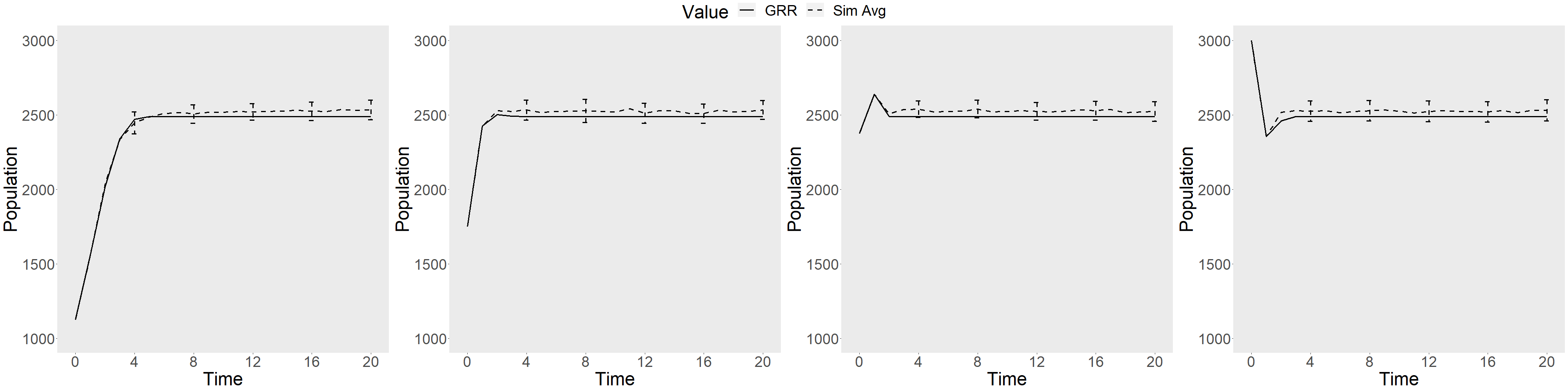}
    \put(0,25){(B)}
    \put(9.5,-2.5){$n_0=1125$}
    \put(34.5,-2.5){$n_0=1750$}
    \put(59.5,-2.5){$n_0=2375$}
    \put(84.5,-2.5){$n_0=3000$}
    \end{overpic}
    \caption{Comparison of GRR population estimates and simulation averages over time for the GoL-like ABM $\cal{A}_\cal{G}=\cal{G}(20,2,8,2,4)$.
    The initial number of living agents ($n_0$) varied uniformly across five values ranging from $500$ to $3000$, altering the initial density of agents across $\Omega_\cal{G}$ for the simulations.
    For all simulations of $\cal{A}_\cal{G}$ considered here, we require that the only agents at time $t=0$ are living agents which start located uniformly at random in the environment $\Omega_\cal{G}$.
    Population averages were generated from 100 simulations for each initial value $n_0$.}
    \label{fig:GoL_AG}
\end{figure}

Note that the ABM $\cal{A}_\cal{G}$ from Section~\ref{ss:ABM} is denoted $\cal{G}(20,$ $2,$ $8,$ $2,$ $4)$ by Definition~\ref{def:GoL}.
Figure~\ref{fig:GoL_AG} shows the first few times of a simulation $\{\cal{X}_t\}_{t\geq 0}$ of $\cal{A}_\cal{G}$; this figure was generated using an implementation\footnote{Our implementation of $\cal{A}_\cal{G}$ in NetLogo expands the environment $\Omega_\cal{G}$ slightly so that it becomes $[0,20]\times[0,20]$.
While this minor modification changes some of the neighborhoods of agents near the boundary of $\Omega_\cal{G}$,
we consider our NetLogo implementations of $\cal{A}_\cal{G}$ and of all the other GoL-like ABMs presented here as being equivalent to their formal descriptions for simplicity.
Of course, one can also modify Definition~\ref{def:GoL} so that it matches the NetLogo implementations.} of $\cal{A}_\cal{G}$ in NetLogo \cite{Wilensky1999}.
Note that the simple rules encoded in $f_\cal{G}$ and $g_\cal{G}$ give rise to clustering despite the fact that surviving agents move randomly within the environment at every time step.

First, observe that the total possible neighborhoods of agents in a GoL-like ABM $\cal{G}(w,$ $\ell_{surv},$ $u_{surv},$ $\ell_{rep},$ $u_{rep})$ partition the environment $\Omega$ into $w^2$ unit squares.

\begin{figure}[ht]
    \centering
    \vspace{\baselineskip}
    
    \begin{overpic}[width=\textwidth]{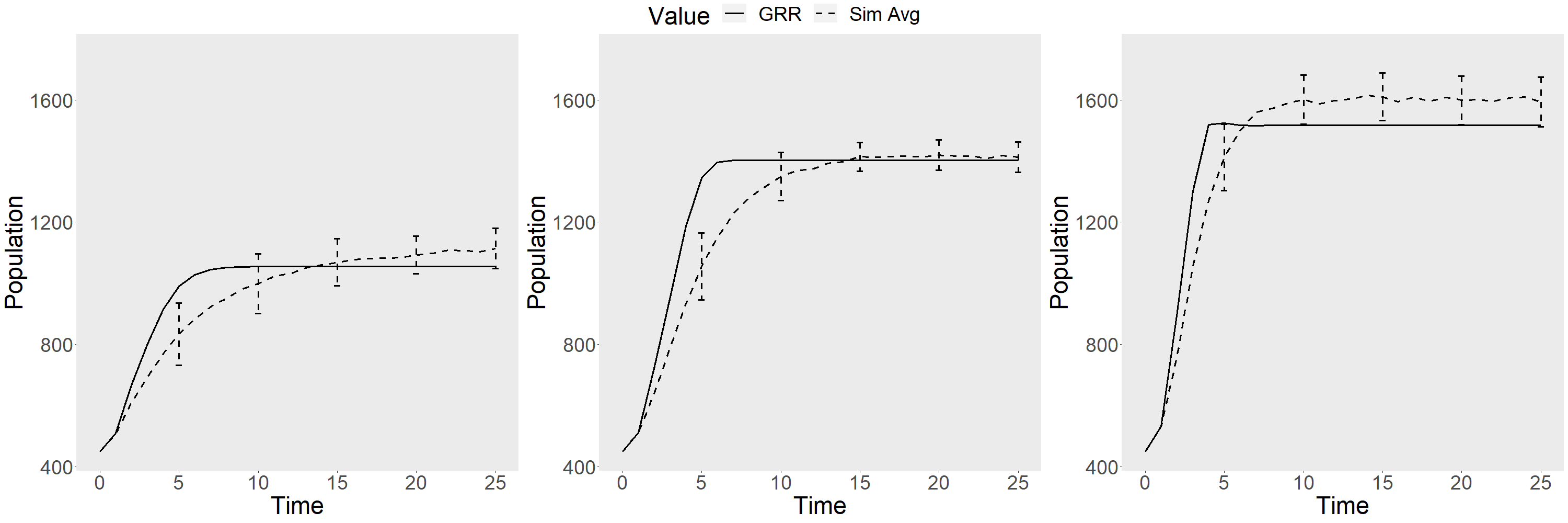}
    \put(0,33){(A)}
    \end{overpic}
    
    \vspace{0.5\baselineskip}
    
    \begin{overpic}[width=\textwidth]{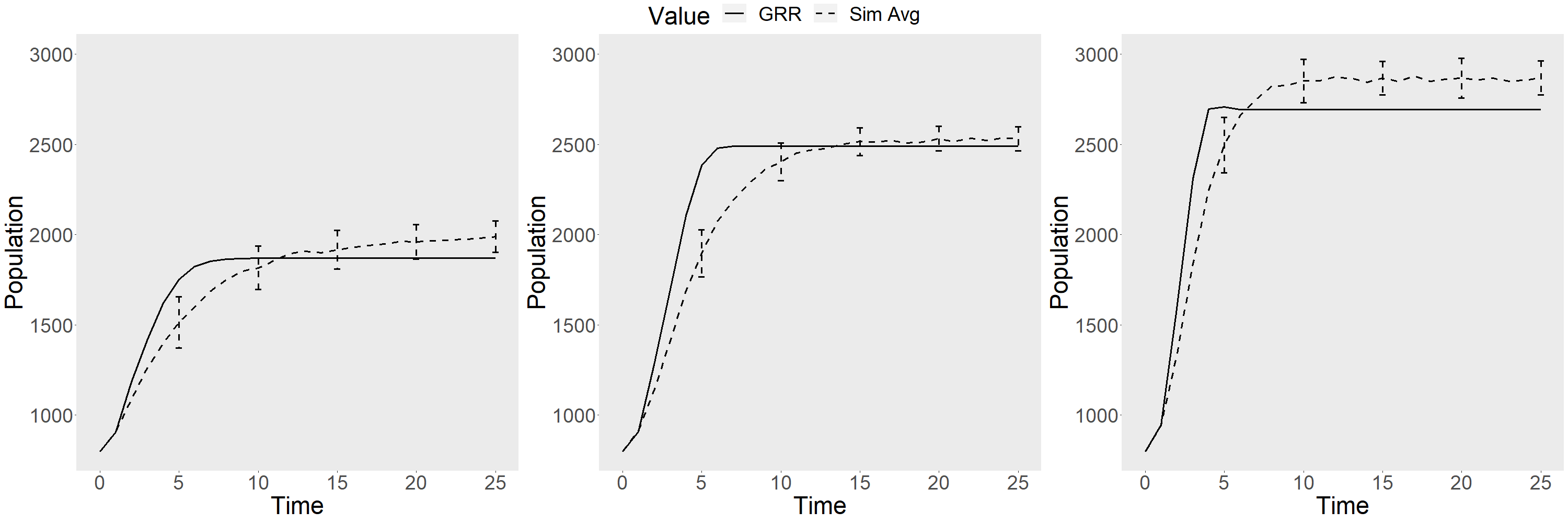}
    \put(0,33){(B)}
    \end{overpic}
    
    \vspace{0.5\baselineskip}
    
    \begin{overpic}[width=\textwidth]{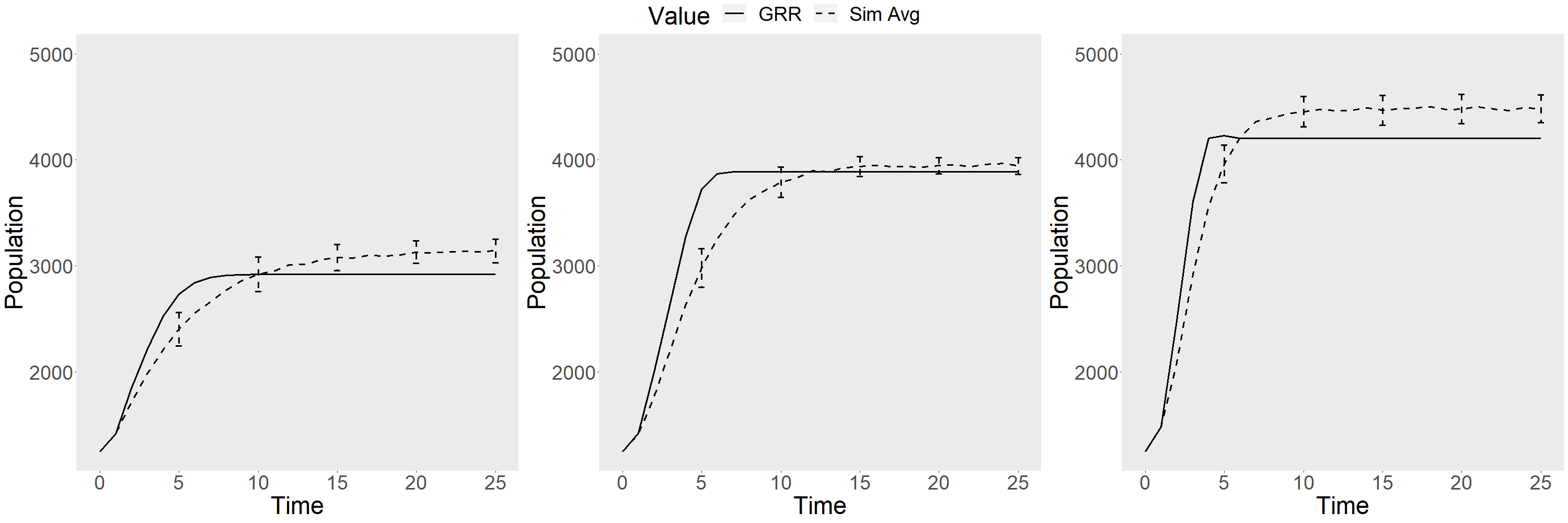}
    \put(0,33){(C)}
    \end{overpic}
    \caption{Comparison of GRR population estimates and simulation averages over time for nine GoL-like ABMs.
    In each ABM $\cal{G}(w,$ $\ell_{surv},$ $u_{surv},$ $\ell_{rep},$ $u_{rep})$, we varied both environment size $w$ and the conditions for survival and reproduction (i.e. $\ell_{surv},$ $u_{surv},$ $\ell_{rep},$ $u_{rep})$).
    The initial number of agents is directly proportional to the size of environment: $n_0=2\cdot w^2$.
    As in Figure~\ref{fig:GoL_AG}, we require that the only agents at time $t=0$ are living agents which start located uniformly at random in their environment.
    Population averages were generated from 100 simulations of each ABM.}
    \label{fig:GoL_variations}
\end{figure}

\section{Rib Development ABM}\label{s:ribABM}
$\quad$

The ABM for early rib development in \cite{Fogel2017}, called ``rib ABM'' $\cal{B}$ hereafter, assessed how genetic modifications regulating development and/or cell proliferation and death affected patterning during rib cage bone formation.
The ribs can be horizontally divided into two compartments: the proximal part connected to the spine and the distal part adjacent to the breastbone.
As the spine and ribs are formed, concentration gradients produced by Hedgehog (Hh) protein diffusing between cells serve as determinants for cells to make their fate decisions between the proximal and distal segments.
The analysis in \cite{Fogel2017} mainly focused on the effects of removing two genes related to the behavior of the agents (cells) called Sonic hedgehog ({\em Shh}) and Apoptotic protease-activating factor 1 ({\em Apaf1}).
These two genes are represented in the agent rules as controlling (i) the rates of cell death and proliferation ({\em Apaf1}) and (ii) the transition to a proximal or distal cell state via the Hh gradient intensity ({\em Shh}). 
In the rib ABM $\cal{B}$, the undetermined (yellow) cells change into proximal (red) or distal (blue) cells depending on the local Hh concentration under four different settings: untreated (``Normal''), {\em Apaf1} knock-out (``{\em Apaf1} KO''), {\em Shh} knock-out (``{\em Shh} KO''), and double knock-out of genes {\em Apaf1} and {\em Shh} (``{\em Apaf1;Shh} DKO'').
This ABM recapitulated the experimental results for each condition, producing different populations of undetermined, proximal, and distal cells.
In this section, we will approximate the GRR of each cell type (yellow, red, and blue, respectively) for $\cal{B}$ under the four aforementioned settings.

\begin{figure}[ht]
\begin{minipage}{.49\textwidth}
\centering
\hspace*{0.06\linewidth}
\frame{\includegraphics[width=0.92\linewidth]{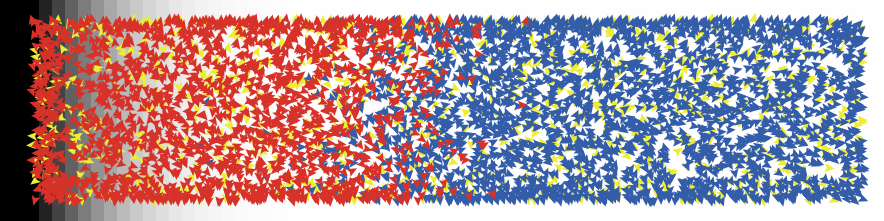}}

\vspace{0.5\baselineskip}
\includegraphics[width=\linewidth]{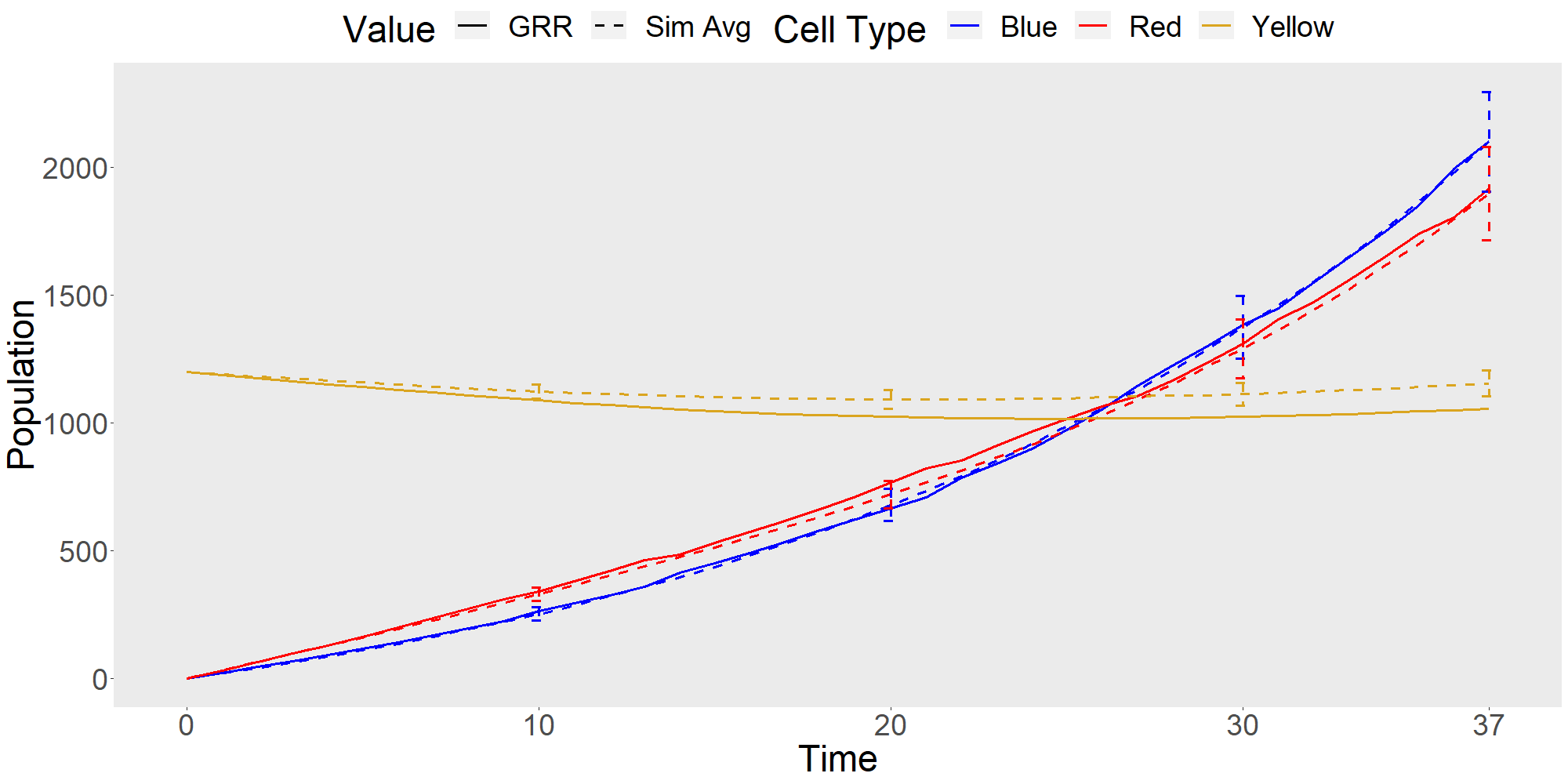}
(a) Normal

\vspace{0.5\baselineskip}

\hspace*{0.06\linewidth}
\frame{\includegraphics[width=0.92\linewidth]{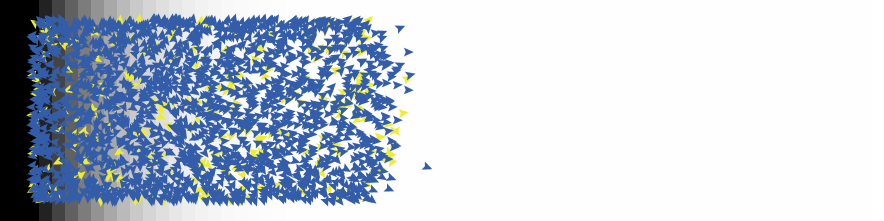}}

\vspace{0.5\baselineskip}
\includegraphics[width=\linewidth]{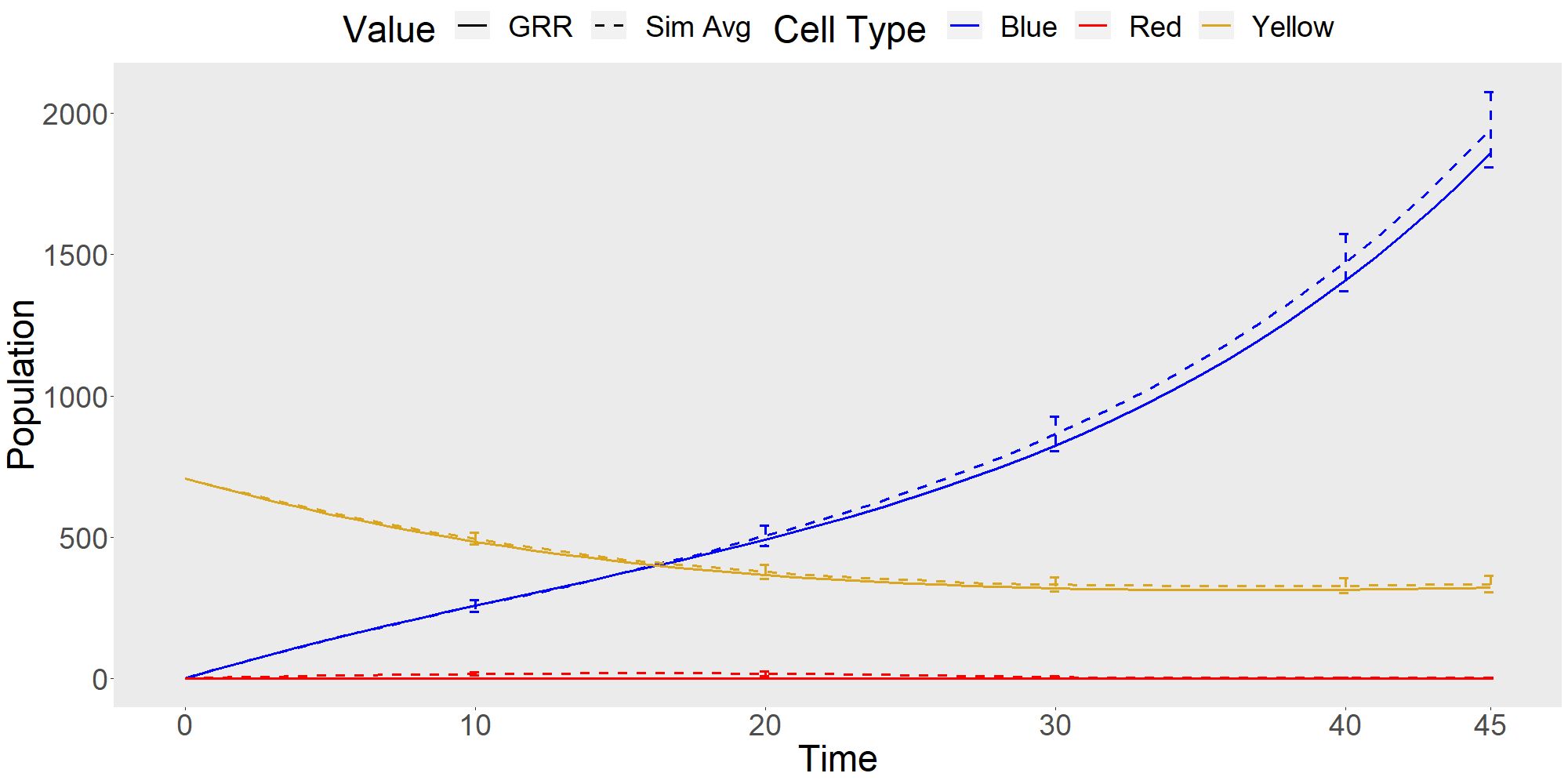}
(c) {\em Shh} KO

\end{minipage}
\hfill
\begin{minipage}{.49\textwidth}
\centering
\hspace*{0.06\linewidth}
\frame{\includegraphics[width=0.92\linewidth]{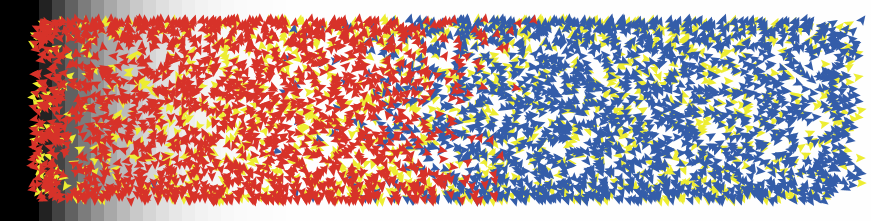}}

\vspace{0.5\baselineskip}
\includegraphics[width=\linewidth]{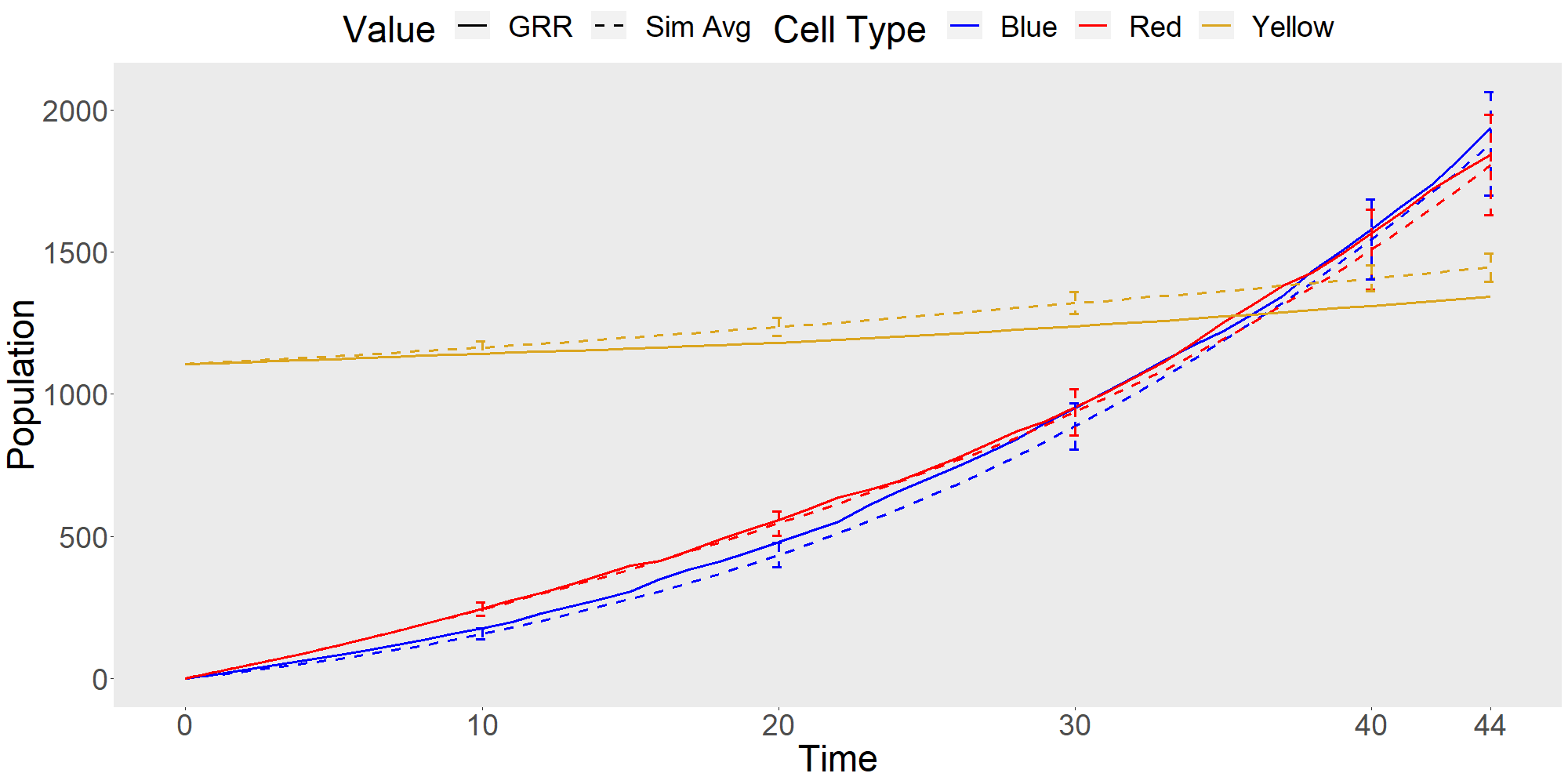}
(b) {\em Apaf1} KO

\vspace{0.5\baselineskip}

\hspace*{0.06\linewidth}
\frame{\includegraphics[width=0.93\linewidth]{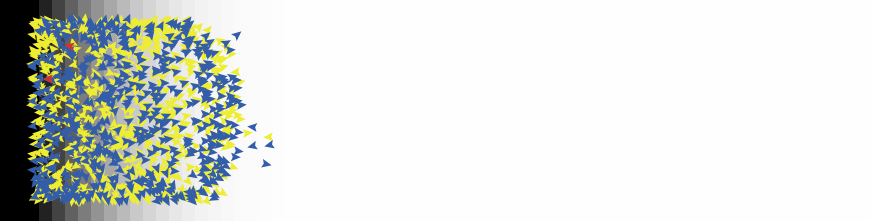}}

\vspace{0.5\baselineskip}
\includegraphics[width=\linewidth]{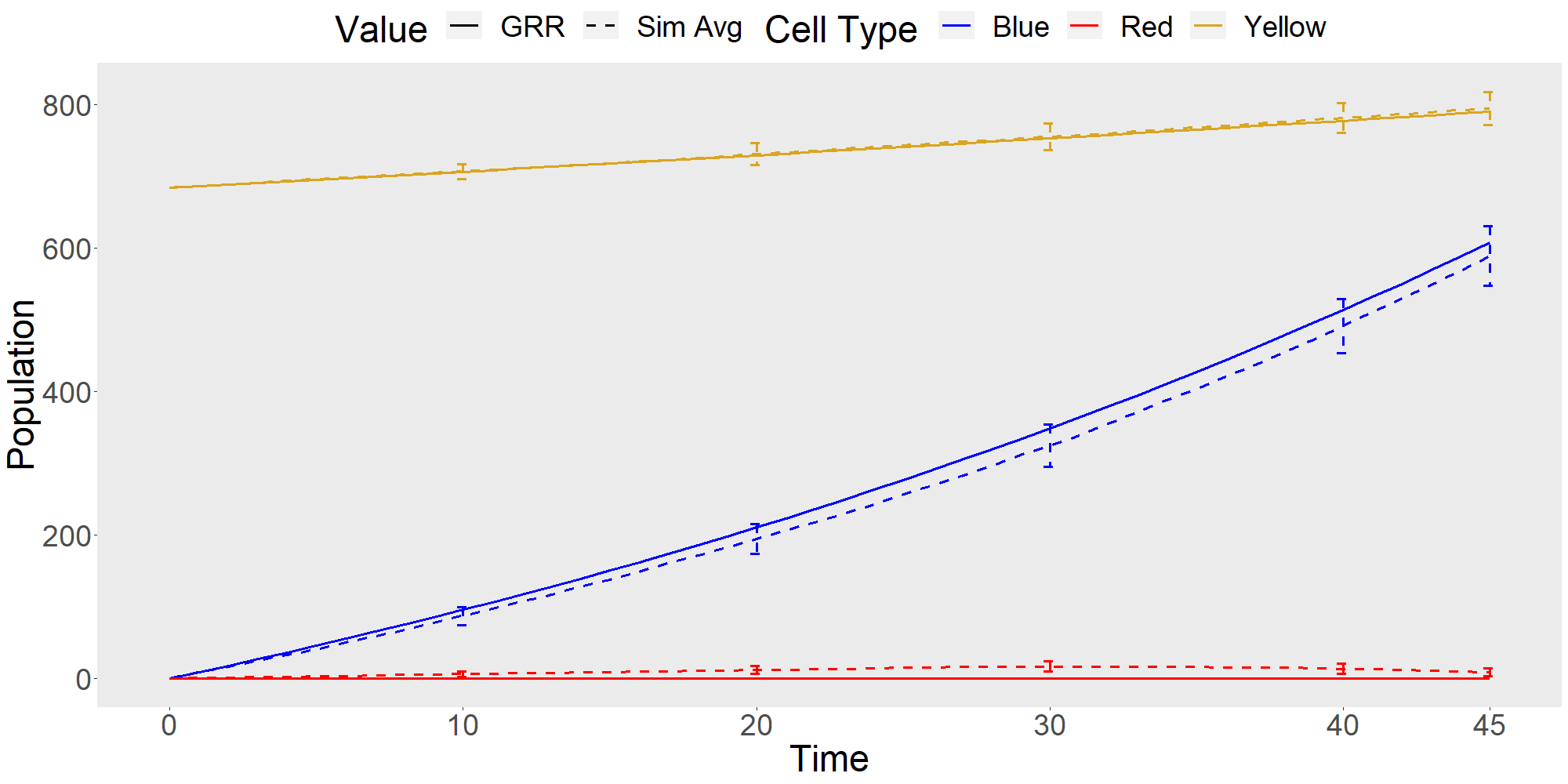}
(d) {\em Apaf1;Shh} DKO

\end{minipage}

\vspace{\baselineskip}
\caption{Comparison of the rib ABM $\cal{B}$ simulations against our GRR estimates for the four settings (i.e. phenotypes) (a) - (d) considered in \cite{Fogel2017}.
For each setting, cell type population averages were generated from 100 simulations.
Each population plot is also visualized with a corresponding simulation snapshot.}
\label{fig:rib_phenotypes}
\end{figure}

\begin{figure}[!ht]
    \begin{overpic}[width=\textwidth]{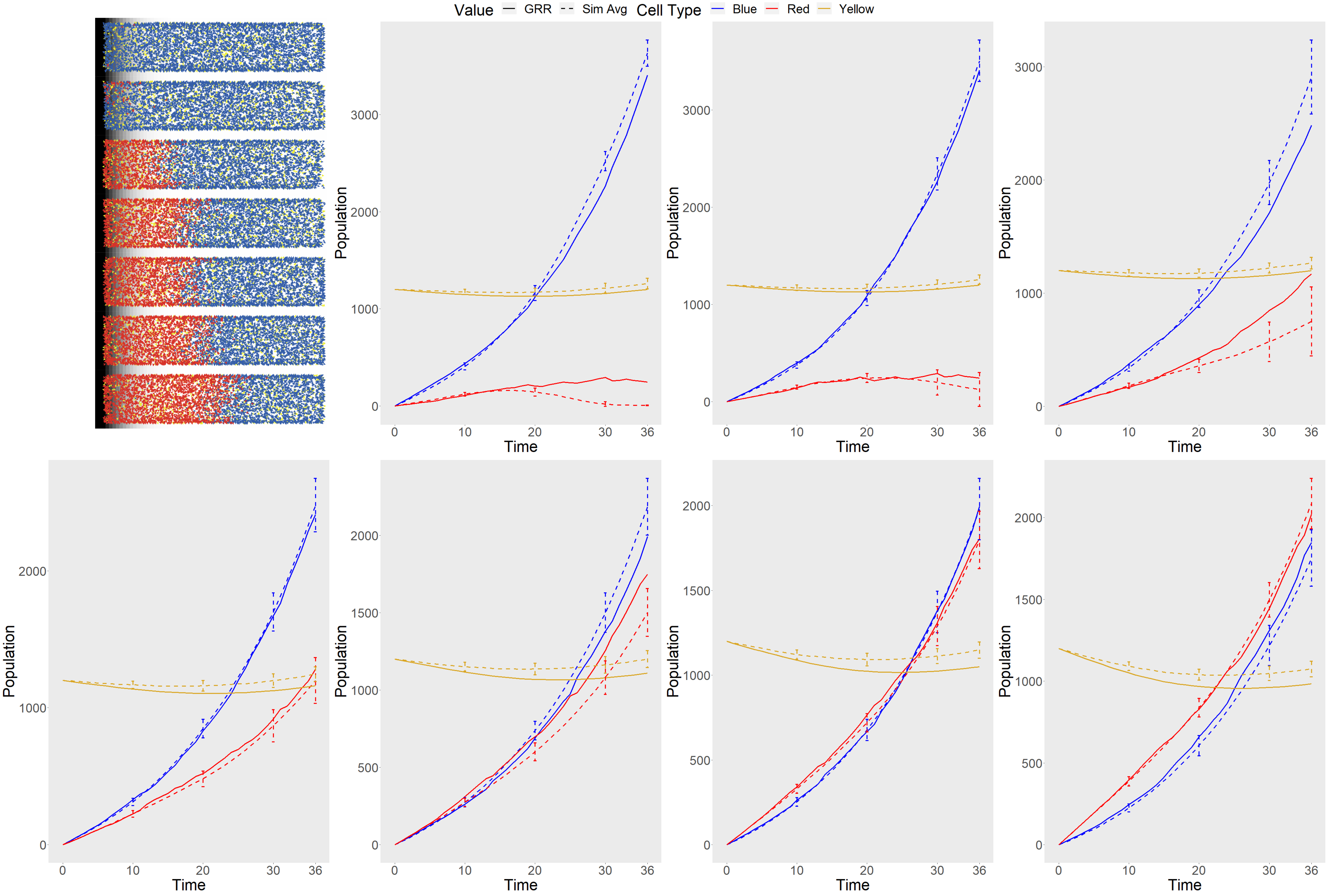}
    \put(1.5, 63){-0.4}
    \put(1.5, 58.5){-0.2}
    \put(4.5, 54.2){0}
    \put(2.5, 49.9){0.2}
    \put(2.5, 45.6){0.4}
    \put(2.5, 41.3){0.6}
    \put(2.5, 37){0.8}
    \put(30,62){\colorbox{white}{-0.4}}
    \put(55,62){\colorbox{white}{-0.2}}
    \put(80,62){\colorbox{white}{0}}
    \put(05,29){\colorbox{white}{0.2}}
    \put(30,29){\colorbox{white}{0.4}}
    \put(55,29){\colorbox{white}{0.6}}
    \put(80,29){\colorbox{white}{0.8}}
    \end{overpic}
    \caption{Comparison of the rib ABM against our GRR estimates for the seven {\em Shh} (log) intensity values considered in \cite{Fogel2017}.
    For each variation, cell type population averages were generated from 100 simulations.
    Each population plot is also visualized with a corresponding simulation snapshot (top left).}
    \label{fig:rib_shh_sweep}
\end{figure}
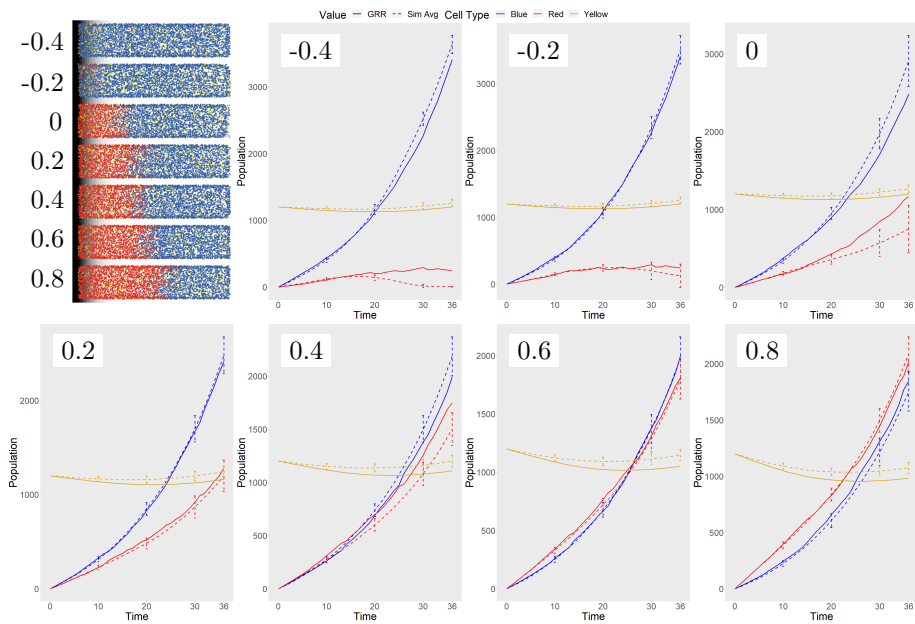

\section{Discussion}\label{s:discussion}
In this work, we considered approximating changes in state densities using our framework for formalizing ABMs. 

\appendix

\section{General Definitions}\label{a:def}
We now generalize the definitions introduced in Section~\ref{s:def} to allow (i) for agents to have shapes (i.e. be represented as connected subsets) and (ii) for local transition and production rules to be stochastic.

\begin{definition}\label{def:gen_agent}
Let $\Sigma$ be a finite set of {\em states}, and let $\Omega$ be (1) a connected, bounded subset of $\bb{R}^n$ or (2) a graph $(V,E)$.
An {\em agent} $a=(s,\cal{P},\cal{N})$ is an ordered triple where $s\in\Sigma$, and $\cal{P},\cal{N}$ are defined as follows with respect to $\Omega$:
\begin{itemize}
    \item $\cal{P}$ is (1) a connected subset of $\Omega$ or (2) a vertex in the graph $\Omega=(V,E)$, as appropriate.
    \item $\cal{N}$ is (1) a connected subset such that $\cal{P}\subseteq\cal{N}\subseteq\Omega$ or (2) $\cal{N}$ is a finite subset such that $p\in\cal{N}\subseteq V$. 
\end{itemize}
In either of the cases (1) or (2), we call $\cal{P}$ the {\em shape of $\alpha$}, $s$ the {\em state of $\alpha$}, and $\cal{N}$ the {\em neighborhood of $\alpha$}.
For brevity, we use the notation $s(\alpha)$ to refer to $s$ and use similar notation for $\cal{P}$ and $\cal{N}$.
Finally, we say that $\Omega$ is an {\em environment} in this context and use $\Lambda(\Sigma,\Omega)$ to denote the set of all possible agents. 
\end{definition}

\bibliographystyle{siamplain}

\end{document}